\documentclass[final]{elsart}

\usepackage{mathrsfs}
\usepackage{latexsym,amsfonts,amssymb}
\usepackage{graphicx,epsfig,subfigure}
\usepackage{psfrag}
\usepackage[sort&compress]{natbib}

\newcommand{\R}{\mathbb R}
\newcommand{\C}{\mathbb C}
\newcommand{\N}{\mathbb N}
\newcommand{\set}[1]{\left\{#1\right\}}
\newcommand{\spec}{\operatorname{{Spec}}}
\newcommand{\trace}{\operatorname{{Trace}}}
\newcommand{\tracep}{\operatorname{{Trace^+}}}

\newcommand{\proof}{{\it Proof}.\,\,}
\newcommand{\eproof}{\hfill $\Box$}

\begin{document}

\begin{frontmatter}

% Title, authors and addresses

% use the thanksref command within \title, \author or \address for footnotes;
% use the corauthref command within \author for corresponding author footnotes;
% use the ead command for the email address,
% and the form \ead[url] for the home page:
% \title{Title\thanksref{label1}}
% \thanks[label1]{}
% \author{Name\corauthref{cor1}\thanksref{label2}}
% \ead{email address}
% \ead[url]{home page}
% \thanks[label2]{}
% \corauth[cor1]{}
% \address{Address\thanksref{label3}}
% \thanks[label3]{}

\title{Asymptotic Stability at Infinity for Differentiable Vector Fields
of the Plane}

% use optional labels to link authors explicitly to addresses:
% \author[label1,label2]{}
% \address[label1]{}
% \address[label2]{}

\author{Carlos Gutierrez\thanksref{aa}},
\thanks[aa]{Partially supported by FAPESP Grant Tem\'atico $\#$ 03/03107-9, and by CNPq Grant $\#$
306992/2003-5.} \ead{gutp@icmc.usp.br}
\author{Benito Pires\thanksref{bb}\corauthref{ca}},
\corauth[ca]{Corresponding author. Fax: +55-16-33739650.}
\thanks[bb]{Supported by FAPESP Grant $\#$ 03/03622-0.}
\ead{bpires@icmc.usp.br}
\author{Roland Rabanal\thanksref{cc}}
\thanks[cc]{Supported by CNPq Grant $\#$ 141853/2001-8.}
\ead{roland@mat.uab.es}

\address{Instituto de Ci\^encias Matem\'aticas e de Computa\c{c}\~{a}o, Universidade de
S\~ao Paulo, Caixa Postal 668, 13560-970, S\~ao Carlos SP, Brazil.}

\begin{abstract}
Let $X:\R^2\backslash\overline{D}_\sigma\to\R^2$ be a differentiable
$($but not necessarily $C^1)$ vector field, where $\sigma>0$ and
$\overline{D}_\sigma=\set{z\in\R^2:\Vert z\Vert\le\sigma}$. Denote
by $\mathcal{R}(z)$ the real part of $z\in\C$. If for some
$\epsilon>0$ and for all $p\in\R^2\backslash\overline{D}_\sigma$, no
eigenvalue of $D_p X$ belongs to $(-\epsilon,0]\cup
\{z\in\C:\mathcal{R}(z)\ge 0\}$, then:\newline
 \noindent a) For all $p\in\R^2\backslash\overline{D}_\sigma$, there
is a unique positive semi--trajectory of $X$ starting at $p$;
\newline\noindent b) It is associated to $X,$  a well defined number
$\mathcal{I}(X)$ of the extended real line $[-\infty,\infty)$
(called  the index of $X$ at infinity) such that for some constant
vector $v\in\R^2$ the following is satisfied: if $\mathcal{I}(X)$ is
less than zero $($resp. greater or equal to zero$)$, then the point
at infinity $\infty$ of the Riemann sphere $\R^2\cup\set{\infty}$ is
a repellor $($resp. an attractor$)$ of the vector field $X+v$.
\end{abstract}

\end{frontmatter}

Submitted to the Journal of Differential Equations on 11 January
2006.

\end{document}